\documentclass[a4paper,10pt,reqno]{amsart}

\usepackage{amsmath,amsfonts,amssymb}
\usepackage{pifont,graphicx,epsfig}
\usepackage[bf]{caption2}
\usepackage{epstopdf}
\usepackage{tikz}

\def\C{\hbox{\font\dubl=msbm10 scaled 1000 {\dubl C}}}
\def\Q{\hbox{\font\dubl=msbm10 scaled 1000 {\dubl Q}}}

\def\Z{\hbox{\font\dubl=msbm10 scaled 1000 {\dubl Z}}}

\def\F{\hbox{\font\dubl=msbm10 scaled 1000 {\dubl F}}}
\def\sF{\hbox{\font\dubl=msbm10 scaled 800 {\dubl F}}}

\title[Descent on Simple Graphs]{A Descent on Simple Graphs\\ -- from Complete to Cycle --\\ and Algebraic Properties of Their Spectra}

\author[Katja M\"onius, J\"orn Steuding, Pascal Stumpf]{Katja M\"onius, J\"orn Steuding, Pascal Stumpf}

\date{January 2018}

\begin{document}

\maketitle

\begin{abstract} 
We investigate a descent on simple graphs, starting with the complete graph on $n$ vertices and ending up with the cycle graph by removing one edge after another. We obtain quantitative results showing that graphs with large diameter must have some eigenvalues of large algebraic degree.
\end{abstract}

{\small \noindent {\sc Keywords:} simple graphs, diameter, eigenvalues, algebraic integers\\
{\sc  Mathematical Subject Classification:} 05C50, 11R04}

\section{Which graphs have integral spectra?}

This question is the title of an interesting article of Harary \& Schwenk from 1973/74; they wrote ``We develop a systematic approach to this question based on operations on graphs. The general problem appears intractable''  \cite{hs}. The spectrum of a graph is the set of eigenvalues of the associated adjacency matrix. The quest of characterizing adjacency matrices for which all eigenvalues are integers seems to be a challenging project, there is no satisfying answer so far. In this paper, we investigate the definitely easier question which graph properties prevent integral eigenvalues. Since all eigenvalues are algebraic integers, it is natural to have a closer look on the algebraic degree of those numbers (i.e., the degree of the minimal polynomial with integer coefficients annihilating this eigenvalue) and ask whether this degree is related to some properties of the underlying graph. It seems that this question has not been studied before.   
\par

For our purpose we investigate a descent on simple graphs with respect to algebraic properties of the eigenvalues of the associated adjacency matrices. It appears that {\it a graph with large diameter has some eigenvalues of large algebraic degree.} This is well illustrated in the table below listing the diameter and the maximum of all algebraic degrees of eigenvalues for the complete graph $K_n$, the complete bipartite graph $K_{m,\ell}$ with $n=m+\ell$, the Payley graph $P(n)$ (existent if $n\equiv 1\bmod\,4$ is a prime power), and the cycle graph $C_n$ on $n$ vertices, respectively. In the following section we shall provide some graph-theoretical background for our study and present the main idea to prove the observation mentioned above.  
$$
\begin{array}{l||c|c}
     & {\rm{diameter}} & {\rm{alg.\,degree}} \\ \hline 
K_n  & 1 & 1 \\
K_{n-\ell,\ell} & \leqslant 2 & \leqslant 2 \\
P(n) & 2 & 2 \\
C_n  & \approx n/2 & \geqslant n/(2\log\log n) 
\end{array}
$$

\section{A descent on graphs}

Let $n\geqslant 3$ be a positive integer. We remove sequentially edges from the complete graph $K_n$ on $n$ vertices until we end up with the cycle graph $C_n$. This leads to a sequence of (undirected) graphs 
$$
G_0:=K_n \to G_1 \to \ldots \to G_m \to G_{m+1} \to \ldots \to G_M:=C_n, 
$$
where $G_{m+1}$ results from $G_m$ by removing one edge and $M=\sharp E(K_n)-\sharp E(C_n)={1\over 2}n(n-3)$. Of course, the sequence of subgraphs $G_m$ of $K_n$ is not uniquely determined. It appears that not every connected undirected graph without loops and multiple edges can be an element of this sequence; for instance, if $n$ is odd, then the sequence cannot contain any $k$-regular graph with odd $k$ since the number of edges is equal to ${1\over 2}kn$. Apart from that there are no further obstructions. McKay \& Wormald \cite{mckay} gave an asymptotic formula for the quite big number of regular graphs. Moreover, for ending up with $C_n$ we may not diminish edges such that the degree of any vertex is less than $2$.  
\bigskip

\begin{center}
  \begin{tikzpicture}
  {
    \def \r {1.0cm};
    \foreach \i in {1,...,5}
    {
      \node [draw, circle, fill = white, inner sep = 0.0cm, minimum size = 0.2cm] 
        (A\i) at ({360 / 5 * (\i - 1)} : \r * 1.0) {};
      \node [draw, circle, fill = white, inner sep = 0.0cm, minimum size = 0.2cm, shift = {(3.1, 0.0)}] 
        (B\i) at ({360 / 5 * (\i - 1)} : \r * 1.0) {};
      \node [draw, circle, fill = white, inner sep = 0.0cm, minimum size = 0.2cm, shift = {(7.2, 0.0)}] 
        (C\i) at ({360 / 5 * (\i - 1)} : \r * 1.0) {};
      \node [draw, circle, fill = white, inner sep = 0.0cm, minimum size = 0.2cm, shift = {(10.3, 0.0)}] 
        (D\i) at ({360 / 5 * (\i - 1)} : \r * 1.0) {};
    }
    \node [fill = white, inner sep = 0.0cm, minimum size = 0.2cm] 
      at (1.6, 0.0) {$\to$};
    \node [fill = white, inner sep = 0.0cm, minimum size = 0.2cm] 
      at (5.2, 0.0) {$\to \dots \to$};
    \node [fill = white, inner sep = 0.0cm, minimum size = 0.2cm] 
      at (8.8, 0.0) {$\to$};
    \path [draw] (A1) -- (A2) -- (A3) -- (A4) -- (A5) -- (A1);
    \path [draw] (A1) -- (A3);
    \path [draw] (A1) -- (A4);
    \path [draw] (A2) -- (A4);
    \path [draw] (A2) -- (A5);
    \path [draw] (A3) -- (A5);
    \path [draw] (B1) -- (B2) -- (B3) -- (B4) -- (B5) -- (B1);
    \path [draw] (B1) -- (B4);
    \path [draw] (B2) -- (B4);
    \path [draw] (B2) -- (B5);
    \path [draw] (B3) -- (B5);
    \path [draw] (C1) -- (C2) -- (C3) -- (C4) -- (C5) -- (C1);
    \path [draw] (C2) -- (C5);
    \path [draw] (D1) -- (D2) -- (D3) -- (D4) -- (D5) -- (D1);
  }
  \end{tikzpicture}
  \\[0.5cm]
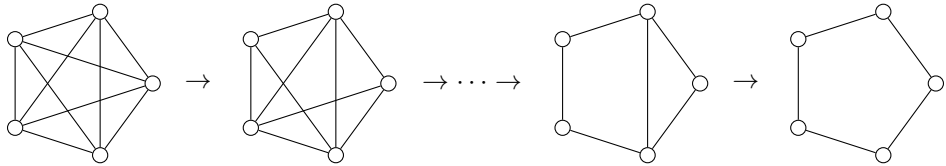
\hspace*{-1.5cm}{\bf Figure 1.} $K_5\ \to\ K_5-\mbox{one edge}\ \to\quad \ldots\quad \to\ C_5+\mbox{one edge}\ \to\ C_5$ 
\end{center}
\medskip

Notice that all eigenvalues are real since the associated adjacency matrices $A(G_m)$ are all symmetric. It follows from the Perron--Frobenius theorem that the largest eigenvalue $\lambda_1^{(m)}$ from the spectrum of the adjacency matrix of $G_m$ is decreasing with increasing $m$, that is (see \cite{brouwer}, p. 33). For $k$-regular graphs it is well known that the largest eigenvalue is equal to $k$, and this eigenvalue is simple if the graph is connected (as in our descent). Moreover, the least eigenvalue is bounded below by $-k$ with equality if, and only if, the graph is bipartite. Since $K_n$ is regular of degree $n-1$ and $C_n$ is regular of degree $2$, 
$$
{\rm{spec}}(C_n)\subset [-2,2]\qquad \mbox{and}\qquad  {\rm{spec}}(K_n)\subset (-(n-1),n-1].
$$
In the general case, if the maximum degree of the vertices of $G_m$ is $k$, then 
$$
{\rm{spec}}(G_m)\subset [-k,k].
$$
Hence, while descending from $K_n$ to $C_n$ the degree of regularity decreases and so the interval for the possible integer eigenvalues shrinks. Furthermore, the diameter of $G_m$ (i.e., the maximum of the minimal lengths of paths connecting vertices of $G_m$) increases from $1$ in the case of the complete graph $K_n$ to $\lfloor n/2\rfloor$ for the cycle $C_n$, where $\lfloor x\rfloor$ denotes the largest integer $\leqslant x$. Obviously, removing edges cannot lower the diameter. Since a graph with diameter $\delta$ has at least $\delta+1$ distinct eigenvalues (see \cite{brouwer}, p. 5), it thus follows that a necessary condition for a graph $G_m$ of maximum vertex degree $k$ from our sequence to have integral spectrum is 
$$
\delta+1\leqslant \sharp([-k,k]\cap \Z)=2k+1.
$$
Thus, {\it whenever $2k<\delta$ with $k$ being the maximum vertex degree of $G_m$, there exists a non-integral eigenvalue within the spectrum.} Actually, there are at least two distinct non-integral eigenvalues as we shall show now.

\section{A little algebra and some graph spectra}

The eigenvalues of $G$ are the roots of the characteristic polynomial $\chi=\det(XI-A)$, where $A$ is the adjacency matrix of $G$. By definition, $\chi$ is a monic polynomial with integer coefficients. By Gauss' lemma (see \cite{mollin}, p. 367), then also all irreducible factors of $\chi$ in $\Q[X]$ are monic and have integral coefficients. Therefore, all eigenvalues are real algebraic integers. It seems that Hoffman \cite{hoffman} was the first to mention this explicitly; moreover, he noticed that with such an algebraic integral eigenvalue $\lambda$ also all its conjugates (that are the other roots of the irreducible polynomial factor of $\chi$) are eigenvalues too, which implies that every eigenvalue has to be a totally real algebraic integer (which rules out numbers as $1/2$ and $\sqrt[3]{4}$ from the spectrum). In particular, every non-integral eigenvalue $\lambda$ has at least one non-integral companion, namely its conjugates and their number equals the algebraic degree of $\lambda$ and its minimal polynomial, respectively. For this and other basic facts from algebraic number theory we refer to \cite{mollin}; the following graph spectra can be found in \cite{brouwer}.    
\medskip

The complete graph $G_0=K_n$ on $n$ vertices has the simple eigenvalue $n-1$ and another eigenvalue $-1$ of multiplicity $n-1$, which we write 
$$
{\rm{spec}}(K_n)=\{n-1^{[1]},-1^{[n-1]}\}
$$
for short. The first deviation from an integral spectrum in our descent are quadratic irrational eigenvalues, already appearing for the graph $G_1$, that is $K_n$ minus one edge. This class of graphs has also been studied by Stark \& Terras; \cite{terras} in context of zeta-functions associated with graphs. Here we find
\begin{equation}\label{minusedge}
{\rm{spec}}(K_n-e)=\{{\textstyle{1\over 2}}(n-3+\sqrt{n^2+2n-7})^{[1]},{\textstyle{1\over 2}}(n-3-\sqrt{n^2+2n-7})^{[1]},0^{[1]},-1^{[n-3]}\};
\end{equation}
notice that the quadratic irrational eigenvalues are conjugate elements in the ring of integers of the number field $\Q(\sqrt{n^2+2n-7})$. Further well-known examples are the bipartite graphs $K_{m,\ell}$ with spectrum
$$
{\rm{spec}}(K_{m,\ell})=\{\sqrt{m\ell}^{[1]},0^{[n-2]},-\sqrt{m\ell}^{[1]}\}
$$
and the ${1\over 2}(q-1)$-regular Paley graphs $P(q)$ for prime powers $q=p^f\equiv 1\bmod\,4$ having spectrum
$$
{\rm{spec}}(P(q))=\{{\textstyle{1\over 2}}(q-1)^{[1]},{\textstyle{1\over 2}}(-1+\sqrt{q})^{[{1\over 2}(q-1)]},{\textstyle{1\over 2}}(-1-\sqrt{q})^{[{1\over 2}(q-1)]}\}.
$$
Here the eigenvalues are real quadratic irrationals as well (except $K_{m,\ell}$ when $m\ell$ is a perfect square). 
\bigskip

\begin{center}
  \begin{tikzpicture}
  {
    \def \r {2.0cm};
    \foreach \i in {1,...,5}
    {
      \node [draw, circle, fill = white, inner sep = 0.0cm, minimum size = 0.2cm] 
        (A\i) at ({360 / 9 * (\i - 3)} : \r * 1.0) {};
    }
    \foreach \i in {6,...,9}
    {
      \node [draw, circle, fill = gray, inner sep = 0.0cm, minimum size = 0.2cm] 
        (A\i) at ({360 / 9 * (\i - 3)} : \r * 1.0) {};
    }
    \foreach \i in {1,...,5}
    {
      \path [draw] (A\i) -- (A6);
      \path [draw] (A\i) -- (A7);
      \path [draw] (A\i) -- (A8);
      \path [draw] (A\i) -- (A9);
    }
  }
  \end{tikzpicture}
  \hspace{0.5cm}
  \begin{tikzpicture}
  {
    \def \r {2.0cm};
    \foreach \i in {1,...,9}
    {
      \node [draw, circle, fill = white, inner sep = 0.0cm, minimum size = 0.2cm] 
        (A\i) at ({360 / 9 * (\i - 1)} : \r * 1.0) {};
    }
    \path [draw] 
      (A1) -- (A2) -- (A3) -- (A4) -- (A5) -- 
      (A6) -- (A7) -- (A8) -- (A9) -- (A1);
    \path [draw] (A1) -- (A5);
    \path [draw] (A1) -- (A6);
    \path [draw] (A4) -- (A8);
    \path [draw] (A4) -- (A9);
    \path [draw] (A7) -- (A2);
    \path [draw] (A7) -- (A3);
    \path [draw] (A9) -- (A2);
    \path [draw] (A3) -- (A5);
    \path [draw] (A6) -- (A8);
  }
  \end{tikzpicture}
  \\[0.5cm]
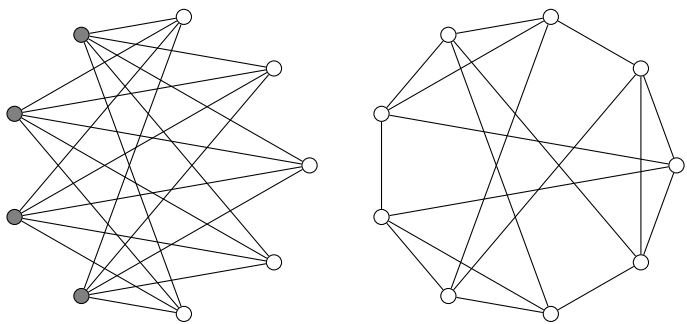
\hspace*{2.7cm}{\bf Figure 2.} $K_{5,4}$ and $P(9)$, both having integral spectrum. 
\end{center}
\medskip

Finally, we end up our descent with the $2$-regular cycle graph $C_n$ on $n$ vertices which has spectrum
\begin{equation}\label{flocke}
{\rm{spec}}(C_n)=\{2\cos(2\pi j/n)\,:\,0\leqslant j<n\};
\end{equation}
all eigenvalues appear with multiplicity 2 except $\lambda =\pm 2$ which are simple, where $\lambda=-2$ is an eigenvalue if, and only if, $C_n$ is bipartite (i.e., when $n$ is even). It appears that all eigenvalues lie in the ring of integers $\Z[\zeta_n]$ of the maximal real subfield $\Q(\zeta_n+\zeta_n^{-1})$ of the cyclotomic field $\Q(\zeta_n)$ of order $n$; here $\zeta_n$ is a primitive $n$th root of unity, $=\exp(2\pi i/n)$ say,  and the degree of the primitive element $\zeta_n+\zeta_n^{-1}$ is ${1\over 2}\varphi(n)$, where $\varphi(n)$ is Euler's totient counting the number of positive integers less than $n$ coprime to $n$. For details about cyclotomy we refer to \cite{wash}.   

\section{Degree inequalities}

Before we continue with counting quadratic irrationals and deducing a criterion for eigenvalues of algebraic degree at least three, we mention another approach to non-integral spectra based on the degrees of the vertices. Given a graph $G$ on $n$ vertices $v_1, \ldots, v_n$, the largest eigenvalue $\lambda_1$ satisfies the following inequalities,
$$
\max_{1\leq j\leq n}\deg (v_j)\geqslant \lambda_1\geqslant {1\over n}\sum_{j\leq n}\deg (v_j)
$$
(see \cite{brouwer}, p. 33); if $G$ is a regular graph minus one edge, $G=K_n-e$ say, it follows from the theorem of Perron-Frobenius that $\lambda_1$ is strictly smaller than $n-1$ (the regularity degree of $K_n$). Moreover,
$$
\sum_{j\leq n}\deg (v_j)=(n-2)(n-1)+2(n-2)=(n-2)(n+1),
$$
hence 
$$
n-2<{(n-2)(n+1)\over n}\leqslant\lambda_1<n-1.
$$
Thus, $\lambda_1$ cannot be integral (as we already know from (\ref{minusedge})). Of course, we do not get any information about the algebraic degree of this eigenvalue. 
\par

Using in addition that the trace of the associated adjacency matrix $A(G)$ is zero (so all eigenvalues add up to zero), that the trace of $A(G)^2$ equals twice the number of edges, and further relations of the trace of powers of $A(G)$ (see \cite{harar}, p. 151, 158), one can deduce that the graph with least number of vertices having at least one eigenvalue with algebraic degree strictly larger than $2$ is the cycle graph $C_5$ on five vertices plus one edge (see Figure 1 for its appearance in our descent). In fact, for $n=5$ all possible spectra consisting of totally real algebraic integers of degree $\leq 5$ are ruled out subsequentially by the conditions on $A(G)^j$ for $0\leq j\leq 4$. It turns out that $C_5$ plus one edge has two integral and three cubic eigenvalues; the corresponding characteristic polynomial is $\chi=X(X+2)(X^3+2X^2-2X-2)$.

\section{Counting quadratic algebraic integers}

For a positive integer $k$ let $a(d,k)$ denote the {\it number of totally real algebraic integers $\alpha$ of degree $d$ for which $\alpha$ and all its conjugates lie inside the interval $[-k,k]$.} In view of $a(1,k)=2k+1$ we derived in Section 1 the existence of an eigenvalue which is an algebraic integer of degree at least $2$ (and according to the degree they come in pairs, triples and so on) whenever the diameter $\delta$ exceeds $2k$. Next we shall compute an upper bound for $a(2,k)$ in order to compute when in our descent the size of the diameter will force the number of distinct eigenvalues to exceed the number of integers and real quadratic irrationals in the spectrum interval.
\par

For this purpose we consider monic quadratic polynomials with integer coefficients, 
$$
P=X^2-bX+c=(X-\alpha_1)(X-\alpha_2),
$$
say. We shall count those polynomials for which the roots $\alpha_1,\alpha_2$ lie inside the interval $[-k,k]$. In view of  
$$
b=\alpha_1+\alpha_2\in[-2k,2k]\ ,\quad c=\alpha_1\alpha_2\in [-k^2,k^2],
$$
the number $a(2,k)$ of totally real quadratic algebraic integers satisfies
$$
a(1,k)+a(2,k)\leqslant 2(4k+1)(2k^2+1)=16k^3+4k^2+8k+2.
$$
In order to count totally real algebraic integers $\alpha_1,\alpha_2$ of degree $2$ the discriminant $b^2-4c$ has to be positive and not a perfect square. Taking this into account one could easily obtain a better bound, for the sake of simplicity, however, we leave this to the interested reader. We observe, that $a(2,k)\leq 17k^3$ for all sufficiently large $k$. 

{\it Assume that $G$ is a graph with diameter $\delta$ and the maximum vertex degree is at most $k$. If $\delta+1>a(1,k)+a(2,k)$, then there are some algebraic integers in the spectrum whose algebraic degrees are at least $3$.} In view of our rough bounds a sufficient inequality for this to happen is $\delta\geq 17k^3$ for large $k$. This inequality is not necessary for an eigenvalue with algebraic degree at least $3$. 

\section{Some algebraic number theory}

For the general case we shall now estimate the number $a(d,k)$ of totally real algebraic integers $\alpha$ of degree $d$ such that $\alpha$ and all its conjugates lie inside the interval $[-k,k]$. Recall that every algebraic integer $\alpha$ of degree $d$ is the root of a monic irreducible polynomial $P_\alpha$ with integer coefficients of degree $d$. Obviously, $P_\alpha$ is uniquely determined and its $d$ roots are called the conjugates of $\alpha$ which we denote by $\alpha_1:=\alpha,\alpha_2,\ldots,\alpha_d$. Thus, 
$$
P_\alpha=\prod_{1\leqslant j\leqslant d}(X-\alpha_j).
$$
The discriminant of $P_\alpha$ is given by
$$
{\rm{disc}}(P_\alpha)=\prod_{i<j}(\alpha_i-\alpha_j)^2.
$$
Since all conjugates lie in the interval $[-k,k]$, we have $\vert \alpha_i-\alpha_j\vert\leqslant 2k$ and
\begin{equation}\label{upper}
{\rm{disc}}(P_\alpha)\leqslant (2k)^{d(d-1)}.
\end{equation}
Hence, {\it if ${\rm{disc}}(P_\alpha)>(2k)^{d(d-1)/2}$, at least one of the conjugates of $\alpha$ lies outside $[-k,k]$.} 
\par

Since a rough bound for the quantity $a(d,k)$ in question is sufficient for our purposes (and, as we shall see below, a tight bound is much harder to find), we drop the condition on $\alpha$ to be totally real. 
\par

Moreover, since every algebraic integer of degree $d$ is an element of a number field $\F$ of dimension (or degree) $d$, considered as a $\Q$-vector space, we may first count all algebraic integers $\alpha$ within a fixed $\F$ of degree $d$ having all conjugates in $[-k,k]$. 
\par

If $\theta_1,\ldots,\theta_d$ denote the $d$ different embeddings of the number field $\F$ into $\C$, then the conjugates of $\alpha$ are given by $\theta_1(\alpha),\ldots,\theta_d(\alpha)$ which, of course, are also algebraic integers of the same degree (but not necessarily elements of $\F$). Let
$$
M=\max\Big\{dk,\left({d\atop 2}\right)k^2,\ldots,\left({d\atop j}\right)k^j,\ldots,k^d\Big\}
$$
and
$$
{\mathcal P}=\Big\{P=X^d+\sum_{0\leqslant j<d}a_jX^j\,:\,a_j\in\Z, \vert a_j\vert \leqslant M\Big\}.
$$
Obviously, ${\mathcal P}$ is finite and consists of $\sharp{\mathcal P}=(2M+1)^d$ elements. Then
$$
{\mathcal A}:=\{\alpha\in \F\,:\,P(\alpha)=0\ \mbox{for some}\ P\in{\mathcal P}\}
$$
is finite too with
$$
\sharp{\mathcal A}\leqslant d\cdot \sharp{\mathcal P}=d\cdot (2M+1)^d.
$$
Denote by $s_1,\ldots,s_d$ the elementary symmetric polynomials of order $d$. If $\alpha\in \F$ satisfies $\vert \theta_j(\alpha)\vert \leqslant k$ for $j=1,\ldots,d$, then
$$
\vert s_j(\theta_1(\alpha),\ldots,\theta_d(\alpha))\vert\leqslant M
$$
for all $j$. Since $\alpha$ is an algebraic integer, all $s_j(\theta_1(\alpha),\ldots,\theta_d(\alpha))$ are rational integers and so are the coefficients of the minimal polynomial of $\alpha$ (since the coefficients are elementary symmetric polynomials of the roots). Thus, $P_\alpha\in{\mathcal P}$ which implies $\alpha\in{\mathcal A}$. (This is essentially Theorem 2.37 in Mollin's monograph \cite{mollin}, p. 90, providing bounds on absolute values.) 
\par 

Since 
$$
M<\sum_{j=0}^d\left({d\atop j}\right)k^j=(k+1)^d,
$$
the {\it number $b(\F,k)$ of algebraic integers $\alpha$ in the field $\F$ such that $\vert \theta_j(\alpha)\vert\leqslant k$ for $j=1,\ldots,d$} satisfies
\begin{equation}\label{a}
b(\F,k)\leqslant d\cdot (2(k+1)^d+1)^d.
\end{equation}
\par

In order to bound $a(d,k)$ it remains to count how many fields $\F$ of degree $d$ make contributions. For this purpose we first notice that every number field has a primitive element, say $\F=\Q(\alpha)$ (so we do not forget to count any $\alpha$), and observe that fields $\F$ with a discriminant too large make no contributions at all. 
\par

The discriminant of a number field $\F$ is equal to the discriminant of any integral basis ${\mathcal B}$ for $\F$, i.e., $\Delta_{\sF}={\rm{disc}}\,({\mathcal B})$. If ${\mathcal B}$ consists of the elements $\alpha_1,\ldots,\alpha_d$, then its discriminant is given by 
$$
{\rm{disc}}\,({\mathcal B})=\det(\theta_j(\alpha_i))^2.
$$
We may form a basis from the consecutive powers of $\alpha$. Then the discriminant of this basis coincides with the one of the minimal polynomial $P_\alpha$ of $\alpha$, 
$$
{\rm{disc}}(\{1,\alpha,\ldots,\alpha^{d-1}\})={\rm{disc}}(P_\alpha)=\prod_{i<j}(\alpha_i-\alpha_j)^2.
$$
Although in general $\{1,\alpha,\ldots,\alpha^{d-1}\}$ is not an integral basis, in view of 
$$
\Delta_{\sF}\cdot [{\mathcal O}_{\sF}:\Z[\alpha]]^2={\rm{disc}}(\{1,\alpha,\ldots,\alpha^{d-1}\}),
$$
where $[{\mathcal O}_{\sF}:\Z[\alpha]]$ stands for the index of the module generated by the powers of $\alpha$ in the ring ${\mathcal O}_{\sF}$ of algebraic integers in $\F$ (which is a positive integer), and (\ref{upper}) we may restrict our counting on those $\F$ satisfying 
\begin{equation}\label{bound}
\Delta_{\sF}\leqslant (2k)^{d(d-1)}.
\end{equation}
With Minkowski's bound from his geometry of numbers one could also provide a lower bound for the discriminant with respect to the degree $d$, however, this would not improve our rough estimate for $a(d,k)$ significantly.
\par

In order to derive an upper bound for $a(d,k)$ we now just need to count how many extensions $\F$ of $\Q$ with fixed degree $d$ and discriminant $\Delta_{\sF}$ satisfying (\ref{bound}) exist. Unfortunately, this seems to be a difficult task! In fact, it has been conjectured that {\it the number $N(d,x)$ of algebraic number fields $\F$ of given degree $d$ and of discriminant $\vert \Delta_{\sF}\vert\leqslant x$} is asymptotically equal to $N(d,x)\sim c_dx$, where $c_d$ is an absolute constant, depending only on $d$. So far this has been confirmed for quadratic and cubic extensions only, i.e., $\Delta_{\sF}\leqslant 3$. Schmidt \cite{schmidt} obtained the estimate 
$$
N(d,x)\ll x^{(d+2)/4}. 
$$
There are notable improvements on this bound: the present best estimate for aribtrary degree is due to Ellenberg \& Venkatesh \cite{ellenberg}; for the sake of simplicity, however, we prefer to use the less complicated and by far more explicit bound due to Schmidt.
\par

In view of (\ref{a}) and (\ref{bound}) we finally arrive at 
\begin{eqnarray*}
a(d,k)&\leqslant & \sum_{\vert \Delta_{\sF}\vert\leqslant (2k)^{d(d-1)/2}}b(\F,k)\\
&\leqslant & N(d,(2k)^{d(d-1)})\cdot d(2(k+1)^d+1)^d\\
&\ll & (2k)^{(d+2)d(d-1)/4}\cdot d(2(k+1)^d+1)^d.
\end{eqnarray*}
For sufficiently large $k$ this quantity is less than $d2^{d^3/3}k^{d^2}$. And for $d=2$ this bound is of order $k^4$ which is worse than our estimate for the quadratic case, so we shall not expect our estimates to be close to the truth. The point here is just that one can derive such a bound. 
\par

Anyway, {\it if the diameter satisfies $\delta+1>a(1,k)+\ldots+a(d,k)$, then there is some algebraic integer in the spectrum whose algebraic degree is at least $d+1$.} In view of our bounds such an eigenvalue exists if $\delta\geq d2^{d^3/3}k^{d^2}$ for large $k$. So our rough conclusion is that a large diameter forces the spectrum to have some eigenvalues of large algebraic degree. 
\smallskip

A final remark about our descent. Kronecker \cite{kronecker} showed that if $\alpha$ is a non-zero totally real algebraic integer which is not of the form $2\cos(\pi r)$ with rational $r$, then the largest absolute value of the conjugates of $\alpha$ is strictly larger than $2$. Robinson \cite{robinson} showed that Kronecker's result is best possible: for every interval of length larger than $4$, there are infinitely many full sets of conjugates of algebraic integers lying in this interval. The difference to our descent is that the algebraic integers appearing in our descent have bounded algebraic degree.

\section{Graphs with many edges and eigenvalues with large algebraic degree}

Taking our observations on graphs with eigenvalues of large algebraic degree into account we may ask how soon can they appear in our descent? In other words: {\it what is the minimal number of edges we have to remove from $K_n$ to have a graph having an eigenvalue of algebraic degree at least $d$?} 
\par

We cannot answer this question but provide an example (which might be optimal in some sense). For this aim we recall that the complement of a simple graph $G=(V,E)$ is defined as the graph $\overline{G}$ having the identical set of vertices and edges exactly there where $G$ has no edges, i.e., if $\overline{G}=(V,\overline{E})$, then $e\in \overline{E}$ if, and only if, $e\not\in E$. If $G$ is a $k$-regular graph on $n$ vertices and its eigenvalues are given by $\lambda_1=k\geqslant\ldots,\geqslant \lambda_n$, then $\overline{G}$ is $n-k-1$ the eigenvalues of $\overline{G}$ are given by
$$
n-k-1\geqslant -1-\lambda_n\geqslant \ldots\geqslant -1-\lambda_2
$$
(for this basic theorem from spectral graph theory see \cite{brouwer}, p. 4).
\bigskip
\begin{center}
  \begin{tikzpicture}
  {
    \def \r {1.2cm};
    \foreach \i in {1,...,5}
    {
      \node [draw, circle, fill = white, inner sep = 0.0cm, minimum size = 0.2cm] 
        (A\i) at ({360 / 5 * (\i - 1)} : \r * 1.0) {};
      \node [draw, circle, fill = white, inner sep = 0.0cm, minimum size = 0.2cm, shift = {(3.5, 0.0)}] 
        (B\i) at ({360 / 5 * (\i - 1)} : \r * 1.0) {};
      \node [draw, circle, fill = white, inner sep = 0.0cm, minimum size = 0.2cm, shift = {(7.0, 0.0)}] 
        (C\i) at ({360 / 5 * (\i - 1)} : \r * 1.0) {};
    }
    \node [fill = white, inner sep = 0.0cm, minimum size = 0.2cm] 
      at (1.8, 0.0) {$=$};
    \node [fill = white, inner sep = 0.0cm, minimum size = 0.2cm] 
      at (5.3, 0.0) {$+$};
    \path [draw] (A1) -- (A2) -- (A3) -- (A4) -- (A5) -- (A1);
    \path [draw] (A1) -- (A3);
    \path [draw] (A1) -- (A4);
    \path [draw] (A2) -- (A4);
    \path [draw] (A2) -- (A5);
    \path [draw] (A3) -- (A5);
    \path [draw] (B1) -- (B2) -- (B3) -- (B4) -- (B5) -- (B1);
    \path [draw] (C1) -- (C3);
    \path [draw] (C1) -- (C4);
    \path [draw] (C2) -- (C4);
    \path [draw] (C2) -- (C5);
    \path [draw] (C3) -- (C5);
  }
  \end{tikzpicture}
  \\[0.5cm]
\hspace*{-4.6cm}{\bf Figure 3.} $K_5\ =\ C_5\ +\ \overline{C_5}$ 
\end{center}
\medskip

If we remove $n$ edges in form of a cycle graph $C_n$ from $K_n$, then we obtain an $n-3$-regular graph $G_n$ with ${1\over 2}n(n-3)$ edges. Since $G_n$ is the complement of $C_n$, it thus follows from (\ref{flocke}) that the spectrum of $G_n$ is given by 
$$
{\rm{spec}}(G_n)={\rm{spec}}(\overline{C_n})=\{n-3,-1-2\cos({2\pi j/n})\,:\,0<j<n\}.
$$
Thus, after removing only $n$ edges one can arrive at a graph having eigenvalues of algebraic degree ${1\over 2}\varphi(n)$. Obviously, we have $\varphi(n)\leq n-1$ with equality if, and only if, $n$ is prime. In the other direction it is known that
$$
\liminf_{n\to\infty}\varphi(n)\cdot {\log\log n\over n}=\exp(-\gamma),
$$
where $\gamma:=\lim_{N\to\infty}(\sum_{n\leq N}{1\over n}-\log N)=0.577\ldots$ is the Euler-Mascheroni constant (see \cite{tenenbaum}, p. 85); this shows that $\varphi(n)$ is {\it almost} of size $n$ (and explains the corresponding entry in the table of the introductory section). Thus, the eigenvalues of a cycle graph $C_n$ have algebraic degree almost as large as possible. {\it Are there graphs having eigenvalues of degree closer to $n$?} 
\medskip

\begin{center}
  \begin{tikzpicture}
  {
    \def \r {2.0cm};
    \foreach \i in {1,...,7}
    {
      \node [draw, circle, fill = white, inner sep = 0.0cm, minimum size = 0.2cm] 
        (A\i) at ({360 / 7 * (\i - 1)} : \r * 1.0) {};
    }
    \path [draw] 
      (A1) -- (A2) -- (A3) -- (A4) -- (A5) -- (A6) -- (A7) -- (A1);
    \path [draw, line width = 1.5pt] (A2) -- (A7);
    \path [draw] (A1) -- (A4);
    \path [draw, line width = 1.5pt] (A1) -- (A5);
  }
  \end{tikzpicture}
  \hspace{0.5cm}
  \begin{tikzpicture}
  {
    \def \r {2.0cm};
    \foreach \i in {1,...,7}
    {
      \node [draw, circle, fill = white, inner sep = 0.0cm, minimum size = 0.2cm] 
        (A\i) at ({360 / 7 * (\i - 1)} : \r * 1.0) {};
    }
    \path [draw] 
      (A1) -- (A2) -- (A3) -- (A4) -- (A5) -- (A6) -- (A7) -- (A1);
    \path [draw, line width = 1.5pt] (A7) -- (A2);
    \path [draw] (A3) -- (A5);
    \path [draw, line width = 1.5pt] (A1) -- (A5);
  }
  \end{tikzpicture}
  \\[0.5cm]
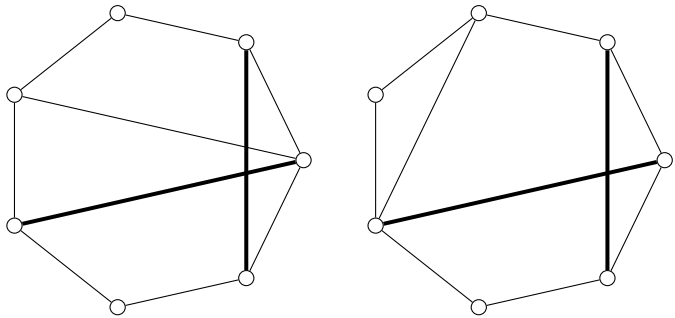
{\bf Figure 4.} Two graphs looking rather similar but having rather different spectra.
\end{center}
\medskip

Figure 4 shows two graphs each of which with $n=7$ vertices and $10$ edges in almost the same configuration. The characteristic polynomial of the graph on the left splits over the rationals completely, so the spectrum is integral: $\{3^{[1]},1^{[2]},0^{[1]},-1^{[1]},-2^{[2]}\}$. On the contrary, the characteristic polynomial $\chi=X^7-10X^5+4X^4+21X^3-8X^2-8X+2$ of the graph on the right is irreducible and therefore all its seven eigenvalues have algebraic degree $7$. 
\par

Can we generalize from this example? Let $\Gamma_n$ be the graph on $n$ vertices with adjacency matrix 
$$
A=A(C_n)+B+B^\top
$$
with $A(C_n)$ being the adjacency matrix of $C_n$ with all entries equal to $1$ next to the diagonal and $B$ is the $n\times n$-matrix with zero entries except for the first row $(0\, \cdots\, 0\, 1\, 1\, 0\, 1\, 0))$; here $B^\top$ denotes the transpose of $B$. Numerical experiments show that the characteristic polynomials $\chi_n$ of $\Gamma_n$ are irreducible for $n=7,8,\ldots,20$, however, we could not find any proof for the case of an arbitrary $n\geq 7$ so far. 

\bigskip

\bigskip

\noindent {\footnotesize Katja M\"onius, J\"orn Steuding, Pascal Stumpf\\ Department of Mathematics, W\"urzburg University, Emil-Fischer-Str. 40, 97\,074 W\"urzburg, Germany, steuding@mathematik.uni-wuerzburg.de}


\begin{thebibliography}{9}

\bibitem{brouwer}{\sc A.E. Brouwer, W.H. Haemers}, {\it Spectra of Graphs}, Springer, 2012

\bibitem{ellenberg}{\sc J.S. Ellenberg, A. Venkatesh}, The number of extensions of a number field with fixed degree and bounded discriminant, {\it Ann. Math.} {\bf 163} (2006), 723-741

\bibitem{harar}{\sc F. Harary}, {\it Graph Theory}, Perseus Books, 1969

\bibitem{hs}{\sc F. Harary, A.J. Schwenk}, Which graph have integral spectra? Graphs and Combin., Proc. Capital Conf., Washington, D. C. 1973, Lect. Notes Math. 406, 45-51 (1974) 

\bibitem{hoffman}{\sc A.J. Hoffman}, Eigenvalues of graphs, in: ``Stud. Graph Theory'', Part II, {\it MAA Stud. Math.} {\bf 12}, 225-245 (1975)

\bibitem{kronecker}{\sc L. Kronecker}, Zwei S\"atze \"uber Gleichungen mit ganzzahligen Coefficienten, {\it J. reine angew. Math.} {\bf 53} (1857), 173-175

\bibitem{mckay}{\sc B.D. McKay, N.C. Wormald}, Asymptotic enumeration by degree sequences of graphs with degrees $o(n^{1/2})$, {\it Combinatorica} {\bf 11} (1991), 369-382

\bibitem{mollin}{\sc R.A. Mollin}, {\it Algebraic number theory}, Boca Raton, Chapman and Hall/CRC 1999

\bibitem{robinson}{\sc R.M. Robinson}, Intervals containing infinitely many sets of conjugate algebraic integers, in: {\it Studies in Mathematical Analysis and Related Topics: Essays in Honor of George P\'olya}, Stanford 1962, 305-315 

\bibitem{schmidt}{\sc W.M. Schmidt}, Number fields of given degree and bounded discriminant, in: {\it Columbia University number theory seminar}, New York, 1992. Paris: Soci\'et\'e Math\'ematique de France, {\it Ast\'erisque} {\bf 228} (1995), 189-195

\bibitem{terras}{\sc H.M. Stark, A.A. Terras}, Zeta functions of finite graphs and coverings. III. {\it Adv. Math.} {\bf 208}, 467-489 (2007) 

\bibitem{tenenbaum}{\sc G. Tenenbaum}, {\it Introduction to Analytic and Probabilistic Number Theory}, Cambridge University Press, 1995

\bibitem{wash}{\sc L.C. Washington}, {\it Introduction to cyclotomic fields}, Springer, 1997, 2nd ed. 

\end{thebibliography}
\end{document}